\documentclass[11pt]{article}

\usepackage{enumerate}
\usepackage{mathrsfs}
\usepackage{amsmath, color}
\usepackage{latexsym}
\usepackage{amssymb}
\usepackage{amsthm}
\usepackage{amsfonts}
\usepackage{dsfont}
\usepackage{hyperref}
\usepackage{xcolor,cancel}

\newtheorem{dfn}{Definition}[section]
\newtheorem{thm}[dfn]{Theorem}
\newtheorem{rmk}[dfn]{Remark}

\newtheorem{prop}[dfn]{Proposition}

\theoremstyle{definition}

\numberwithin{equation}{section}
\newcommand{\IE}{{\mathbb{E}}}
\newcommand{\IP}{{\mathbb{P}}}
\newcommand{\IR}{{\mathbb{R}}}
\newcommand{\FF}{{\mathcal{F}}}
\newcommand{\EE}{{\mathcal{E}}}

\def\wt{\widetilde}
\def\<{\langle}
\def\>{\rangle}

\title{\bf  On Transition Density Functions of Skew Brownian Motions with Two-Valued Drift}
\date{\today}
\author{{\bf Shuwen Lou}}

\begin{document}

\maketitle

\begin{abstract}
In this article, based on the results in \cite{GS},  we derive the explicit full family of transition density functions of skew Brownian motion (SBM in abbreviation)  with two-valued drift for all $t>0$.  As an important step of this result, it is also shown in this paper that SBM with two-valued drift is a strong Markov process by finding its symmetrizing measure and canonical scale function, from which one can tell what values of the drift make such a process transient or recurrent. 

\end{abstract}

\medskip
\noindent
{\bf AMS 2010 Mathematics Subject Classification}:  60J35,  60J45, 60J65, 60J60.

\smallskip\noindent
{\bf Keywords and phrases}: Skew Brownian motion,  two-valued drift,   Dirichlet forms,  $h$-transform, transition density functions.

\section{Introduction}\label{Intro}

Let $(\Omega, \FF, \IP)$ be a probability space and let $\{B_t, \FF_t, t\ge 0\}$  be a standard Brownian motion (BM in abbreviation) with its natural filtration. Consider the following SDE:
\begin{equation}\label{sbm-2valued-drift}
dY_t = dB_t +m(Y_t) dt +(2p-1) dL_t^0, \quad Y_0=y_0,
\end{equation}
where  $y_0\in \IR$, $0<p<1$,
\begin{equation}\label{def-m(x)}
m(x)= m_1\cdot \mathbf{1}_{\{x\ge 0\}}+m_2\cdot \mathbf{1}_{\{x<0\}},  \quad m_1, m_2\in \IR, 
\end{equation}
and $(L^0_t)_{t\ge 0}$ is the symmetric semimartingale local time of $Y$  at $0$  in the sense that 
\begin{equation*}
L^0_t =\lim_{\epsilon \downarrow 0}\frac{1}{2\epsilon}\int_{0}^t  \mathbf{1}_{\{Y_s\in [-\epsilon, \epsilon]\}} ds.
\end{equation*}  
Such an SDE is known to have a unique weak  solution  (see the first paragraph of the proof to Theorem \ref{thm-potential-theory-sbm-2-valued} for details), which is called an {\it SBM with a two-valued drift}. In the classic work \cite{KS1}, standard Brownian motion with a two-valued drift  (i.e., $p=1/2$ in \eqref{sbm-2valued-drift})  was studied, where both the transition density functions and the joint distributions of the process itself and its functionals were obtained. Similar results were obtained for SBM with a constant drift (i.e., $m_1=m_2$ in \eqref{def-m(x)})  in \cite{ABTWW, ABTWW1}.

In financial mathematics, SBM with a two-valued drift can be viewed as the rescaled process of options following two-valued  local volatility model (two-valued LVM), see, e.g., \cite{GS} and the references therein. In practice, option prices under LVM can be numerically calculated by solving numerically the corresponding
partial differential equations. Because of this, finding the explicit transition density functions of $Y$ which can be numerically approximated is of great interest.   In \cite[Theorem 2]{GS}, using the method of discrete approximation, the    explicit density function of the position  of   SBM  with two-valued drift starting from $0$  is obtained. To state their result, denote by $p(t,0, x)$ the density function of the position of $Y$ starting from $0$ with respect to the Lebesgue measure, i.e., for any Borel measurable set $A\subset \IR$ and any $t>0$,
\begin{equation*}
\int_A p(t, 0, x)dx =\IP^0\left[Y_t\in A\right],
\end{equation*}
where $\IP^0$ stands for the distribution of $Y$ starting from $Y_0=0$. 

As a continuation of \cite{GS}, in this paper we further answer a few questions unanswered in \cite{GS}. It is proved in this paper that SBM with two-valued drifts are continuous strong Markov processes, and we present their symmetrizing measures and canonical scale functions  in Theorem \ref{thm-potential-theory-sbm-2-valued}.  Based upon this,  as the main result of this paper, we derive the full family of transition density functions $\{p(t,x,y)\}_{t>0, x,y\in \IR}$ of $Y$ as the following theorem states.

\begin{thm}\label{main-result}
Let $m_*:=|m_1|\wedge |m_2|$, and $m^*=|m_1|\vee |m_2|$. Denote the transition densities of  $Y$ with respect to Lebesgue measure by $p(t,x,y)$.   It holds 
\begin{enumerate}
\item[\rm (i)]  For $x> 0,\, y\le 0$, 
\begin{align*}
&\frac{2(1-p)}{\sqrt{2\pi t^3}} e^{-m_1x+m_2y-\frac{(m^*)^2t}{2}  }\int_0^\infty (l-y+x) e^{-lpm_1+l(1-p)m_2 -\frac{(l-y+x)^2}{2t}} dl\le p(t,x,y) 
\\
&\le  \frac{2(1-p)}{\sqrt{2\pi t^3}} e^{-m_1x+m_2y-\frac{m_*^2t}{2} }   \int_0^\infty (l-y+x) e^{-lpm_1+l(1-p)m_2 -\frac{(l-y+x)^2}{2t}}dl;
\end{align*}

\item[\rm (ii)]  For $x> 0,\, y>0$, 
\begin{align*}
&\frac{2p}{\sqrt{2\pi t^3}}e^{m_1(y-x)-\frac{(m^*)^2t}{2}}\int_0^\infty (l+y+x)e^{-lm_1p+lm_2(1-p)-\frac{(l+y+x)^2}{2t}}dl
\\
&+\frac{e^{-\frac{m_1^2}{2}+m_1(y-x)}}{\sqrt{2\pi t}}\left( e^{-\frac{|x-y|^2}{2t}} - e^{-\frac{|x+y|^2}{2t}}  \right)\le p(t,x,y)
\\
& \le \frac{2p}{\sqrt{2\pi t^3}}e^{m_1(y-x)-\frac{m_*^2t}{2}}\int_0^\infty (l+y+x)e^{-lm_1p+lm_2(1-p)-\frac{(l+y+x)^2}{2t}}dl
\\
&+\frac{e^{-\frac{m_1^2}{2}+m_1(y-x)}}{\sqrt{2\pi t}}\left( e^{-\frac{|x-y|^2}{2t}} - e^{-\frac{|x+y|^2}{2t}}  \right);
\end{align*}

\item[\rm (iii)] For $x<0,\, y\ge 0$,
\begin{align*}
&   \frac{2p}{\sqrt{2\pi t^3 }}  e^{-m_2x+m_1y-\frac{(m^*)^2t}{2}} \int_0^\infty (l+y-x) e^{-lpm_1+l(1-p)m_2-\frac{(l+y-x)^2}{2t}}dl \le p(t,x,y) 
\\
&\le \frac{2p}{\sqrt{2\pi t^3 }} e^{-m_2x+m_1y-\frac{m_*^2t}{2}} \int_0^\infty (l+y-x) e^{-lpm_1+l(1-p)m_2-\frac{(l+y-x)^2}{2t}} dl;
\end{align*}

\item[\rm (iv)] For $x<0,\, y< 0$,
\begin{align*}
&\frac{2(1-p)}{\sqrt{2\pi t^3}}e^{m_2(y-x)-\frac{(m^*)^2t}{2}} \int_{0}^\infty (l-y-x)e^{-lm_1p+lm_2(1-p)-\frac{(l-y-x)^2}{2t}}dl
\\
 &+\frac{e^{-\frac{m_2^2}{2}+m_2(y-x)}}{\sqrt{2\pi t}}\left( e^{-\frac{|x-y|^2}{2t}} - e^{-\frac{|x+y|^2}{2t}}  \right)\le p(t,x,y)
 \\
& \le \frac{2(1-p)}{\sqrt{2\pi t^3}}e^{m_2(y-x)-\frac{m_*^2t}{2}}\int_{0}^\infty (l-y-x)e^{-lm_1p+lm_2(1-p)-\frac{(l-y-x)^2}{2t}}dl
\\
&+\frac{e^{-\frac{m_2^2}{2}+m_2(y-x)}}{\sqrt{2\pi t}}\left( e^{-\frac{|x-y|^2}{2t}} - e^{-\frac{|x+y|^2}{2t}}  \right).
\end{align*}

\end{enumerate}
\end{thm}

\begin{rmk}
\begin{description}
\item[\rm (i)] For each of  the four cases above, when $m^*=m_*$, i.e., $|m_1|=|m_2|$, the upper bound and lower bound match. 
\item[\rm (ii)] Since the case $x=0$ is already covered by \cite[Theorem 2]{GS} (restated in \eqref{density-Y-start-from-zero}),  the entire family of $\{p(t,x,y)_{t>0, x,y\in \IR}\}$ is found.
\item[\rm (iii)] Although $\{p(t,x,y)_{t>0, x,y\in \IR}\}$  can also be expressed explicitly using the Duhamel's principle (see, e.g., \cite[Lemma 3.17]{Lou1} and the proof therein), as $Y$ is just a SBM with an explicit drift, the expression via Duhamel's principle is in the format of an infinite sum of multiple integrals, and the $n^{\rm th}$ integral in the sum is an $n$-fold integral, for $n\ge 1$. From the practical perspective, our result is significantly better, especially since such types of processes arise from financial mathematics. 
\end{description}

\end{rmk}

The rest of this paper is organized as follows. In Section \ref{S:2} we first derive the Dirichlet form and the symmetrizing measure associated with $Y$, which are given in Thereom \ref{thm-potential-theory-sbm-2-valued}, which also states the criteria for the transience or  the  recurrence of $Y$. Proposition \ref{density-existence}    establishes the existence of the transition density function of $Y$ using the classic argument in \cite{QZ}. Proposition \ref{density-of-q} is an important step  to  the proof of our main result, Theorem \ref{main-result}, which is the explicit transition density function of $Y$ killed upon hitting $0$. With all the preparation in Section \ref{S:2}, the proof to  Theorem \ref{main-result} is given in Section \ref{pf-of-main-result}.

In this paper we follow the convention that, for a strong Markov process $X$   with state space $E$, for any compact set $K\subset E$, we define $\sigma_K:=\inf\{t>0, X_t\in K\}$. For any open domain $D\subset E$, we define $\tau_{D}:=\inf\{t>0:  X_t \notin D\}$.   Also as  a convention, for a strong Markov process, we use $\IP^x$ and $\IE^x$ to represent the probabilities and the expectations related to this process starting from $x$ in its state space.

\section{Some Potential Theory for Skew Brownian Motion with Two-Valued Drift}\label{S:2}

Given the SDE in \eqref{sbm-2valued-drift} with the pair $m_1, m_2\in \IR$, in this paper we set 
\begin{equation}
\ell(dx):=e^{2m_1x}dx|_{[0,\infty)} + \frac{(1-p)}{p}e^{2m_2x} dx|_{(-\infty, 0)}.
\end{equation}
The following theorem shows that $\ell$ is the symmetrizing measure for $Y$ by giving its associated Dirichlet form on $L^2(\IR, \ell)$.
\begin{thm}\label{thm-potential-theory-sbm-2-valued}
The SDE given by \eqref{sbm-2valued-drift} is well-posed (i.e., having weak solutions with pathwise uniqueness), and its solution $Y$ is a symmetric diffusion process on $\IR$   associated with the following Dirichlet form $(\EE^Y, \FF^Y)$    on $L^2(\IR, \ell)$: 
\begin{equation}\label{Dirichlet-form-Y}
\begin{aligned}
&	\FF^Y=\{f\in L^2(\IR, \ell):\, f\text{ absolutely continuous on }\IR,  \, f'\in L^2(\IR, \ell)\}, \\
&	\EE^Y(f,g)=\frac{1}{2}\int_\IR f'(x)g'(x)\ell(dx),\quad f,g\in \FF^Y,
\end{aligned}
\end{equation}
where $f'$ stands for the weak derivative  with respect to $\ell$ for all  $f\in \FF^Y$.   Furthermore,  $Y$ is recurrent if and only if both $m_1\le 0$ and $m_2\ge 0$ are satisfied. Otherwise $Y$ is transient. 
\end{thm}
\begin{proof}
Well-posedness: Using the notations in  \cite[Theorem 7.1]{Li}, the SDE given by \eqref{sbm-2valued-drift}  corresponds to 
\begin{equation*}
\rho(x)=e^{2m_1x}\mathbf{1}_{[0,\infty)}(x) + \frac{(1-p)}{p}e^{2m_2x} \mathbf{1}_{(-\infty, 0)}(x)
\end{equation*}
and
\begin{equation*}
\mu(dx)=m_1\cdot dz|_{(0, +\infty)}+m_2\cdot dz|_{(-\infty, 0)}+(2p-1)\delta_{\{0\}}(dz). 
\end{equation*}
It is therefore easy to verify that the conditions in \cite[Theoreom 7.1]{Li} are satisfied. Thus \eqref{sbm-2valued-drift} is well-posed with a unique weak solution.

On account of \cite[Lemma 2.2.7]{CF}, one can tell  that the  bilinear form  $(\EE^Y,\FF^Y)$ given in \eqref{Dirichlet-form-Y} is indeed    a Dirichlet form on $L^2(\IR, \ell)$ with canonical scale function 
\begin{equation*}
s(x)=\int_0^x \left(\frac{\ell (dr)}{dr}\right)^{-1}dr,
\end{equation*}
where $\frac{\ell(dr)}{dr}$ stands for the Radon-Nikodym derivative of $\ell$ with respect to the one-dimensional Lebesgue measure. That is,
\begin{equation}\label{scale-function}
s(x)=\left\{
\begin{aligned}
&-\frac{(1-p)}{p}\int_x^0 e^{-2m_2r}dr, &x<0;
\\
&0, & x=0;
\\
&\int_0^x e^{-2m_1r}dr, & x>0.
\end{aligned}
\right.
\end{equation}
From the above we see that:
\begin{itemize}
\item $s(+\infty)=+\infty$ if and only if $m_1\le 0$. When $m_1< 0$,  for any $c\in (-\infty, +\infty)$,  $\ell((c,+\infty))$ is a positive constant in $(0, +\infty)$. 
\item $s(-\infty)=-\infty$ if and only if $m_2\ge 0$. When $m_2> 0$,  for any $c\in (-\infty, +\infty)$,  $\ell((-\infty, c))$ is a positive constant in $(-\infty, 0)$. 
\end{itemize}
In view of the discussion in \cite[\S2.2.3]{CF}, it follows that 
\begin{itemize}
\item When $m_1\le 0$, $+\infty$ is an non-approachable endpoint for $\IR$. When $m_1>0$, it is an approachable but non-reguar endpoint. 
\item  When $m_2\ge  0$, $-\infty$ is an non-approachable endpoint for $\IR$. When $m_2<0$, it is an approachable but non-reguar endpoint. 
\end{itemize}
It now follows from \cite[Proposition 2.2.10]{CF} that $(\EE^Y, \FF^Y)$ is a regular, strongly local, and irreducible Dirichlet form on $\IR$,  which therefore has an $\ell$-symmetric diffusion process associated with it.  Furthermore, such a process is recurrent if and only if both $m_1\le 0$ and $m_2\ge 0$ are satisfied. Otherwise, it is transient. 

Finally, to complete the proof,  a  standard argument using Fukushima decomposition can show that such a diffusion process is indeed the solution to \eqref{sbm-2valued-drift} in view of its well-posedness. Below we spell out the details for readers' convenience. 
Take $f(r):=r\in \FF^Y_\mathrm{loc}$ and consider the Fukushima's decomposition for $f$:
\[
	f(Y_t)-f(Y_0)=M^f_t+N^f_t. 
\]
The martingale part $M^f$ is determined by its energy measure $\mu_{\<f\>}$ and for any $g\in C_c^\infty(\IR)$ (see \cite[Theorem~5.5.2]{FOT}),
\[
	\int gd\mu_{\<f\>}= 2\EE^Y(fg,f)-\EE^Y(f^2,g)=\int gd\ell. 
\]
It follows that $\mu_{\<f\>}=\ell$ and hence $M^f$ is equivalent to a standard Brownian motion. For the zero-energy part $N^f$, by one-dimensional integration by parts,  we have
\begin{eqnarray}
-\EE^Y(f,g) &= & -\frac{1}{2}\int_{\IR} g'(x)\ell(dx) \nonumber
\\
& = & \frac{(2p-1)}{2p} \cdot g(0)
+m_1\int_0^\infty g(x)e^{2m_1x}dx +\frac{m_2(1-p)}{p}\int_{-\infty}^0 g(x)e^{2m_2x}dx \nonumber
\\
&=& \frac{(2p-1)}{2p} \cdot g(0) + \int_{\IR} g(x)\left( m_1\cdot \mathbf{1}_{\{x>0\}}+ m_2\cdot \mathbf{1}_{\{x<0\}}\right)\ell (dx).
\end{eqnarray}
Thus \cite[Corollary~5.5.1]{FOT} yields that $N^f$ is of bounded variation, and its associated signed smooth measure is
\[
	\mu_{N^f}=\frac{(2p-1)}{2p}\cdot \delta_0+m(x)\ell(dx),
\]
where $m(x)$ is given in \eqref{def-m(x)}. This yields that 
\begin{equation}\label{EQ4YTY}
	Y_t-Y_0=B_t+\int_0^t m(Y_s)ds+\frac{(2p-1)}{2p}\cdot l^0_t,\quad t\geq 0, 
\end{equation}
where $(B_t)$ is a certain standard Brownian motion and $l^0:=(l^0_t)_{t\geq 0}$ is the the local time of $Y$ at $0$, i.e. is the positive continuous additive functional of $M$ having Revuz measure $\delta_0$. Finally, in order to  convert $l^0$ to the corresponding semimartingale local time,  using \cite[Lemma 4.3]{Li} we get 
\begin{equation*}
L_t^0 =\frac{1}{2}\left( 1+\frac{1-p}{p}\right)l^0_t=\frac{1}{2p}l^0_t, 
\end{equation*}
which finishes the proof. 
\end{proof}

\begin{prop}\label{density-existence}
The one-dimensional diffusion process $Y$ admits a family of transition density functions $\{p(t,x,y)\}_{t>0, x,y\in \IR}$ with respect to the Lebesgue measure on $\IR$.
\end{prop}
\begin{proof}
 The diffusion process $Y$ can be obtained from the skew Brownian motion $Z$ characterized by the following SDE  through a drift perturbation
(i.e. Girsanov transform):
\begin{equation*}
dZ_t = dB_t + (2p-1)d L^0_t (Z) ,
\end{equation*}
whose explicit transition density is known, see, e.g., \cite{RY}. In view of the definition of Kato-class functions (see, e.g., \cite[Definition 1.3]{Lou1}), the function $m$ given  in \eqref{def-m(x)} satisfies that  $|m|^2 \in \mathbf{K}_{1}$ on $\IR$. Therefore, the argument in \cite{QZ} shows the existence of the transition densities of $Y$.
\end{proof}
To continue, we first establish the existence and the explicit expression of the  transition density functions of $Y$ killed upon hitting $0$.  In other words,  in the next proposition we find a family of functions $\{q(t,x,y)\}_{t>0, x,y\in \IR\backslash \{0\} }$ such that  for any non-negative function $f\geq 0$ on $\IR\backslash \{0\}$, 
\begin{equation}\label{density-Y-cond-on-hitting-0}
 \int_{\IR\backslash \{0\}}  q (t, x, y)  f(y) dy = \IE^x \left( f(Y_t); t< \sigma_{\{0\}} \right).
\end{equation}
For $x=0$ or $y=0$, we make the convention $q(t,x,y):=0$. The following proposition gives a neat computation for the explicit expression of $q(t,x,y)$ via  an $h$-transform using Dirichlet form theory. For readers' convenience, we record the definition of $h$-transform as follows, which can be found in \cite[Definition 5.5]{Gyrya}.

\begin{dfn}[$h$-transform]
Let $(\mathfrak{E}, \mathfrak{D(E)})$ on $L^2(M, \mu)$ with associated semigroup $(P_t)_{t>0}$ and infinitesimal generator $(L, \mathfrak{D}(L))$. Let $h$ be a measurable positive function on $M$. Let $H$ denote the multiplication by $h$ viewed as a unitary operator
\begin{equation*}
H: L^2(M, h^2d\mu) \rightarrow  L^2(M, d\mu), \quad f\mapsto hf. 
\end{equation*}
 Define $(\mathfrak{E}_h, \mathfrak{D(E_{\mathit{h}})})$, $L_h$, and $P_{h, t}$ to be, respectively, the pulled back closed bilinear form, operator, and semigroup on $L^2(M, h^2d\mu)$
\begin{eqnarray}
\mathfrak{E}_h (f, g) &=& \mathfrak{E} (hf, hg), \quad  \mathfrak{D(E_{\mathit{h}})}=H^{-1}\mathfrak{D(E)};
\\
L_h &=& H^{-1}\circ L \circ H, \quad \mathfrak{D}(L_h) =H^{-1}\mathfrak{D}(L);
\\
P_{h, t} &=& H^{-1}\circ P_t \circ H.
\end{eqnarray}
$(\mathfrak{E}_h, \mathfrak{D(E_{\mathit{h}})})$ is called the h-transform of  $(\mathfrak{E}, \mathfrak{D(E)})$.
\end{dfn}
From this definition one can see that when $(\mathfrak{E}_h, \mathfrak{D(E_{\mathit{h}})})$ is Markovian (i.e., normal contractions operate on it), it is a Dirichlet form on $L^2(M, h^2d\mu)$ associated with self-adjoint semigroup $(P_{h, t})_{t>0}$ and self-adjoint infinitesimal generator $L_h$. 

\begin{prop}
Let $\{q(t,x,y)\}_{t>0, x,y\in \IR\backslash \{0\}}$ be as defined in \eqref{density-Y-cond-on-hitting-0}. 
\begin{equation}\label{density-of-q}
q(t,x,y)=\frac{e^{-\frac{m_1^2}{2}-m_1x+m_1y}}{\sqrt{2\pi t}}\left( e^{-\frac{|x-y|^2}{2t}} - e^{-\frac{|x+y|^2}{2t}}  \right), \quad x,y\in \IR_+, t>0
\end{equation}
\end{prop}

\begin{proof}
In this proof, we let $\IR_+$ be $(0, +\infty)$ not including $0$.
We denote by $Y^{\IR_+}$ the part  process  of $Y$ on $\IR_+$, i.e., killed upon hitting $(-\infty, 0]$. $Y^{\IR_+}$ is characterized by the following Dirichlet form  $(\FF^Y_{\IR_+}, \EE^Y_{\IR_+})$   on $L^2(\IR_+, \ell)$:
\begin{equation}
\begin{aligned}
&\FF^{Y}_{\IR_+}=\{f\in L^2(\IR, \ell):\, f\text{ absolutely continuous on }\IR,\, f'\in L^2(\IR, \ell),\, f=0 \, \,\text{q.e. on }(-\infty, 0]\}, \\
& \EE^{Y}_{\IR_+}(f,g)=\frac{1}{2}\int_{\IR_+} f'(x)g'(x)\ell(dx),\quad f,g\in \FF^{Y}_{\IR_+}.
\end{aligned}
\end{equation}
Therefore, by setting $\phi (x):=e^{m_1x}$ on $\IR_+$, the Dirichlet form $(\FF^{Y}_{\IR_+}, \EE^{Y}_{\IR_+} )$ can be viewed as the $h$-transform of the following Dirichlet form $(\mathcal{G}, \mathcal{A})$ on $L^2(\IR_+, dx)$   with  $h=\phi$: 
\begin{equation}
\begin{aligned}
	\mathcal{G}&=\{u\in L^2(\IR, dx):\, u/\phi \in \FF^{Y}_{\IR_+}\}, \\
	\mathcal{A}(u, v)&=\frac{1}{2}\int_{\IR_+} \left(\frac{u}{\phi}\right)' \left(    \frac{v}{\phi}\right)' \phi ^2(x)dx,\quad u, v\in \mathcal{G},
\end{aligned}
\end{equation}
Denote  the Dirichlet form associated with the standard Brownian motion killed upon exiting $\IR_+$ by $(W^{1,2}(\IR_+),  \frac{1}{2}\mathbf{D})$. It is known that  $C_c^\infty (\IR_+)$ is a core of $( W^{1,2}(\IR_+), \frac{1}{2}\mathbf{D} )$. Since $\phi$ is  positive and $C^\infty$ on $\IR_+$, we know $C_c^\infty (\IR_+)$ is also a core of $(\mathcal{G}, \mathcal{A})$.  Now take any $u, v\in C^\infty_c(\IR_+)$, 
\begin{eqnarray}
\mathcal{A}(u, v) &=& \frac{1}{2}\int_{\IR_+}\frac{u'\phi -u\phi'}{\phi^2}\cdot  \frac{v'\phi -v\phi'}{\phi^2}\cdot \phi^2 dx \nonumber
\\
&=& \frac{1}{2} \int_{\IR_+} \left( u'(x)e^{m_1x}-u(x)\cdot m_1e^{m_1x}  \right)\left( v'(x)e^{m_1x}-v(x)\cdot m_1e^{m_1x} \right)e^{-2m_1x}dx  \nonumber
\\
&=&\frac{1}{2}\int_{\IR_+} u'(x)v'(x)dx + \frac{m_1^2}{2} \int_{\IR_+} u(x)v(x)dx
-\frac{m_1}{2}\int_{\IR_+}\left( u(x)v'(x)+u(x)v'(x) \right) dx \nonumber
\\
&=& \frac{1}{2}\int_{\IR_+} u'(x)v'(x)dx + \frac{m_1^2}{2} \int_{\IR_+} u(x)v(x)dx.
\end{eqnarray}
This implies that the Dirichlet form $(\mathcal{G}, \mathcal{A})$ corresponds to part Brownian motion on $\IR_+$ subject to a constant killing rate $m_1^2/2$. In view of the transition density of part Brownian motion on $\IR_+$ (see, e.g., \cite[Exercise 1.15]{RY}), we know the transition density functions of $(\mathcal{G}, \mathcal{A})$ with respect to the Lebesgue measure are
\begin{equation}
p^{(\mathcal{G}, \mathcal{A})}(t,x,y)=\frac{e^{-\frac{m_1^2}{2}}}{\sqrt{2\pi t}}\left( e^{-\frac{|x-y|^2}{2t}} - e^{-\frac{|x+y|^2}{2t}}  \right), \quad x,y\in \IR_+, t>0.
\end{equation}
Thus by \cite[Lemma 5.6]{Gyrya},  the transition density functions of $Y^{\IR_+}$ with respect to $\ell$   are
\begin{equation*}
p^{Y^{\IR_+}}(t,x,y)=\frac{p^{(\mathcal{G}, \mathcal{A})}(t,x,y)}{\phi (x)\phi (y)}=\frac{e^{-\frac{m_1^2}{2}-m_1(x+y)}}{\sqrt{2\pi t}}\left( e^{-\frac{|x-y|^2}{2t}} - e^{-\frac{|x+y|^2}{2t}}  \right), \quad x,y\in \IR_+, t>0.
\end{equation*}
Finally, since $q(t,x,y)$ stands for the transition density functions of $Y^{\IR_+}$ with respect to Lebesgue measure, the proof is complete by a change of measure.
\end{proof}

\section{Proof of Theorem 1.1}\label{pf-of-main-result}

First of all, in view of the fact that starting from $x>0$, before hitting $0$, $Y$ has the same distribution as a one-dimensional Brownian motion with a constant drift $m_1$ on $\IR$, it follows immediately from \cite[p.197, (5.12)]{KS} that  starting from $x>0$, the hitting time $\sigma_{\{0\}}:= \inf \{t>0: Y_t =0\}$  has the following distribution density:
\begin{equation}\label{hitting-time-dist}
\IP^x\left[ \sigma_{\{0\}}\in ds \right]=\frac{|x|}{\sqrt{2\pi s^3}} e^{-\frac{(m_1s+x)^2}{2s}}, \quad x\in \IR_+, s>0.
\end{equation}
By the definition of $q(t,x,y)$ in \eqref{density-Y-cond-on-hitting-0}, setting 
\begin{equation}\label{decomp-p(t,x,y)}
\overline{p}(t,x,y)=p(t,x,y)-q(t,x,y),
\end{equation}
it holds for any non-negative function $f\geq 0$ on $\IR$ that
\begin{equation}\label{def-overline-p}
 \int_{\IR} \overline p (t, x, y)  f(y) dy = \IE^x \left[ f(Y_t); t\geq \sigma_{\{0\}} \right],
\end{equation}
i.e., $\overline{p}(t,x,\cdot )$ is the density of $Y$ at time $t$ starting from $x$,  provided that the path hits  $0$ before $t$. Thus by the Markov property of $Y$, 
\begin{equation}\label{markov-overline-p}
\overline{p}(t,x,y)=\int_0^t p(t-s, 0, y)\IP^x\left[ \sigma_{\{0\}} \in ds \right].
\end{equation}
Let $p(t,0,x,l)$ be the joint density of $(Y_t, L^{0}_t)$ given $Y_0=0$.    \cite[Theorem 2]{GS} gives that
\begin{align}\label{density-Y-start-from-zero}
&p(t,0,x,l)  = \int_{0}^t \int_{0}^u 
  \frac{a(x) \cdot  l^2 \cdot p (1-p)\cdot |x|}{\sqrt{2}\left( \pi v (u-v)(t-u)\right)^{3/2}}   \exp\bigg\{-\frac{l^2}{2}\left(\frac{p^2}{v}+\frac{(1-p)^2}{u-v} \right) \nonumber
\\
&-\frac{x^2}{2(t-u)}-\frac{m_1^2v+m_2^2(t-v)}{2}-lp(m_1+m_2)+lm_2 +x\cdot m(x) \bigg\}  dvdu, \, t>0, x\in \IR,\, l>0,
\end{align}
where $a(x):= p\mathbf{1}_{\{x> 0\}}+(1-p)\mathbf{1}_{\{x\le 0\}}$. Denote by  $Y^0$  the solution to \eqref{sbm-2valued-drift} when $m_1=m_2=0$,  i.e., 
\begin{equation}
dY^0_t = dB_t  +(2p-1) dL_t^{0,0}, \quad Y^0_0=0,
\end{equation}
where $L^{0,0}_t$ stands for its symmetric local time at $0$. Let $p^0(t,0,x,l)$ be the joint density of $(Y^0_t, L^{0,0}_t)$ given $Y^0_0=0$. Also let $\sigma^0_{\{0\}}:= \inf \{t>0: Y^0_t =0\}$. In the following, by comparing $p(t,0,x,l)$ with $p^0(t,0,x,l)$,  we relate $p(t,x,y)$ to $p^0(t,x,y, l)$, which has explicit expressions given by \cite[Corollary 3.3]{ABTWW, ABTWW1}.    We first claim that:
\begin{equation}\label{claim-local-time-strong-additive}
\IP^x\left[Y_t\in dy, L_t^0\in dl, t> \sigma_{\{0\}}\right]=\int_{s=0}^t \IP^0\left[Y_{t-s}\in dy,  L^0_{t-s}\in dl\right]\IP^x\left[\sigma_{\{0\}}\in ds\right].
\end{equation}
 To justify \eqref{claim-local-time-strong-additive}, we first note by \cite[(1.12)]{BG} and the fact that local times are perfect additive functions (or \cite[Chapter X, Proposition 1.2]{RY}, in which the claim is made for Brownian local times, but the same proof can be adapted to our case) that for any positive random variable $S$ and any $x\in \IR$,
\begin{equation*}
L^0_{S+\sigma_{\{0\}}}=L^0_{\sigma_{\{0\}}}+L^0_S(\theta_{\sigma_{\{0\}}})=L^0_S(\theta_{\sigma_{\{0\}}}),\quad \IP^x\text{-a.s.}
\end{equation*}
Thus by the strong Markov property of $Y$, 
\begin{eqnarray}
&& \IP^x\left[Y_t\in dy, L_t^0\in dl, t> \sigma_{\{0\}}\right] \nonumber
\\
&=& \IP^x\left[ \IE\left[Y_t \in dy, L_{t-\sigma_{\{0\}}}\circ \theta_{\sigma_{\{0\}}} \in dl\Big|\FF_{\sigma_{\{0\}}}\right]; t>\sigma_{\{0\}}    \right] \nonumber
\\
&=&\IP^x\left[ \IP^0\left[Y_{t-\sigma_{\{0\}}}\in dy, L_{t-\sigma_{\{0\}}}\in dl  \right]; t>\sigma_{\{0\}} \right] \nonumber
\\
&=& \int_{s=0}^t \IP^0\left[Y_{t-s}\in dy,  L^0_{t-s}\in dl\right]\IP^x\left[\sigma_{\{0\}}\in ds\right].  \nonumber
\end{eqnarray}
This justifies \eqref{claim-local-time-strong-additive}. In the following we divide our discussion into four cases, depending on whether $x,y$ are greater than $0$ or not. 
\\
{\it Case 1.} $x>0, y\le 0$.  For this case, on account of the continuity and Markov property of $Y$,
\begin{equation}\label{Markov}
p(t,x,y)=\int_0^t p(t-s, 0, y)\IP^x\left[ \sigma_{\{0\}} \in ds \right]. 
\end{equation}
In \eqref{density-Y-start-from-zero}, noticing that $\frac{m_1^2v+m_2^2(t-v)}{2}\ge \frac{m_*^2t}{2}$, we have for all $x\in \IR$,
\begin{eqnarray}
&& p(t,0,x,l) \nonumber
\\
& \le & e^{-\frac{m_*^2t}{2}-lp(m_1+m_2)+lm_2+x\cdot m(x)}  \int_0^t \int_0^u \frac{a(x)   l^2  p (1-p) |x|}{\sqrt{2}\left( \pi v (u-v)(t-u)\right)^{3/2}} e^{-\frac{l^2}{2}\left(\frac{p^2}{v}+\frac{(1-p)^2}{u-v}\right)-\frac{x^2}{2(t-u)} } dvdu \nonumber
\\
& = & e^{-\frac{m_*^2t}{2}-lp(m_1+m_2)+lm_2+x\cdot m(x)}\cdot   p^0(t,0,x,l).\label{p(t,0,x,l)-upper-bound}
\end{eqnarray}
Similarly,
\begin{equation}\label{p(t,0,x,l)-lower-bound}
p(t,0,x,l)\ge e^{-\frac{(m^*)^2t}{2}-lp(m_1+m_2)+lm_2+x\cdot m(x)}\cdot   p^0(t,0,x,l).
\end{equation}
Therefore, the right hand side of \eqref{Markov} can be rewritten as
\begin{eqnarray}
p(t,x,y) 
& = & \int_0^t p(t-s, 0, y)\IP^x\left[ \sigma_{\{0\}} \in ds \right]\nonumber
\\
&\stackrel{\text{Fubini}}{=} &\int_0^\infty\int_0^t p(t-s, 0, y, l)\IP^x\left[ \sigma_{\{0\}} \in ds \right]dl \nonumber
\\
&\stackrel{\eqref{p(t,0,x,l)-upper-bound}}{\le }& \int_0^\infty \int_0^t e^{-\frac{m_*^2(t-s)}{2}-lp(m_1+m_2)+lm_2+y\cdot m(y)} \cdot p^0(t-s,0,x,l) \IP_x[\sigma_{\{0\}}\in ds]dl \nonumber
\\
&\stackrel{\eqref{hitting-time-dist}}{=} &  \int_0^\infty \int_0^t e^{-\frac{m_*^2(t-s)}{2}-lp(m_1+m_2)+lm_2+m_2y}\cdot  p^0(t-s,0,y,l)   \frac{|x|}{\sqrt{2\pi s^3}} e^{-\frac{(m_1s+x)^2}{2s}}  dsdl \nonumber
\\
& = & \int_0^\infty \int_0^t e^{-\frac{m_*^2(t-s)}{2}-lp(m_1+m_2)+lm_2+m_2y}  \cdot p^0(t-s,0,y,l)   \frac{|x|}{\sqrt{2\pi s^3}} e^{-\frac{m_1^2 s}{2}-m_1x-\frac{x^2}{2s}}  ds dl \nonumber
\\
& \le & \int_0^\infty e^{-\frac{m_*^2t}{2}-lp(m_1+m_2)+lm_2+m_2y-m_1x} \int_0^t  p^0(t-s, 0, y, l) \frac{|x|}{\sqrt{2\pi s^3}}e^{-\frac{x^2}{2s}}ds dl \nonumber
\\
& = & \int_0^\infty  e^{-\frac{m_*^2t}{2}-lp(m_1+m_2)+lm_2+m_2y-m_1x}\int_0^t  p^0(t-s,0,y,l) \IP_x\left[\sigma^0_{\{0\}}\in ds\right] dl\nonumber
\\
&\stackrel{\eqref{claim-local-time-strong-additive}}{=}& \int_0^\infty  e^{-\frac{m_*^2t}{2}-lp(m_1+m_2)+lm_2+m_2y-m_1x} p^0(t,x,y,l) dl.
\label{eq3.5}
\end{eqnarray}
In view of \eqref{p(t,0,x,l)-lower-bound}, the lower bound for $p(t,x,y,l)$ can be shown similarly by replacing all the $m_*$  in \eqref{eq3.5}   with $m^*$, i.e., $p(t,x,y)$ has the following two-sided bounds:
\begin{eqnarray}
&&\int_0^\infty  e^{-\frac{(m^*)^2t}{2}-lp(m_1+m_2)+lm_2+m_2y-m_1x} p^0(t,x,y,l) dl \nonumber
\\
&\le & p(t,x,y) \le \int_0^\infty e^{-\frac{m_*^2t}{2}-lp(m_1+m_2)+lm_2+m_2y-m_1x} p^0(t,x,y,l)dl.
\end{eqnarray}
The conclusion  for this case follows  in view of \cite[Corollary 3.3]{ABTWW, ABTWW1} which gives
\begin{equation*}
p^0(t,x,y,l)=\frac{2(1-p)(l-y+x)}{\sqrt{2\pi t^3}} e^{-\frac{(l-y+x)^2}{2t}}, \quad l>0, y\le 0.
\end{equation*}
{\it Case 2.} $x>0, y>0$. For this case,
\begin{align}\label{case2-density}
\IP^x[Y_t\in dy]&=\IP[Y_t\in dy, \sigma_{\{0\}}\ge t]+\IP^x[Y_t\in dy, \sigma_{\{0\}}< t]\nonumber
\\ 
&\stackrel{\eqref{density-Y-cond-on-hitting-0}}{=}q(t,x,y)dy +\IP^x [Y_t\in dy, \sigma_{\{0\}}< t] \nonumber
\\
& \stackrel{\eqref{markov-overline-p}}{=} q(t,x,y)dy+\int_0^t  \IP^0[Y_{t-s}\in dy]\; \IP^x[\sigma_{\{0\}} \in ds]. 
\end{align}
The first term on the right hand side above can be replaced with \eqref{density-of-q}. For the second term, similar to the computation for \eqref{eq3.5},
\begin{eqnarray}
&&  \int_0^t  \IP^0[Y_{t-s}\in dy]\; \IP^x[\sigma_{\{0\}} \in ds]  \nonumber
\\
& = & \int_0^t \int_0^\infty p(t-s, 0, y, l)dl dy \; \IP^x\left[ \sigma_{\{0\}} \in ds \right] \nonumber
\\
&\stackrel{\eqref{p(t,0,x,l)-upper-bound}}{\le } &  \int_0^t \int_0^\infty e^{-\frac{m_*^2(t-s)}{2}-lp(m_1+m_2)+lm_2+m_1y}\cdot   p^0(t-s,0,x,l)  dldy\; \IP_x[\sigma_{\{0\}}\in ds] \nonumber
\\
%&= & \int_0^\infty \int_0^t   e^{-\frac{m_*^2(t-s)}{2}-lp(m_1+m_2)+lm_2+m_1 y}\cdot   p^0(t-s,0,y,l)   \frac{|x|}{\sqrt{2\pi s^3}} e^{-\frac{(m_1s+x)^2}{2s}}  dsdl  dy\nonumber
%\\
%& = & \int_0^t \int_0^\infty e^{-\frac{m_*^2(t-s)}{2}-lp(m_1+m_2)+lm_2+m_2y}\cdot   p^0(t-s,0,y,l)   \frac{|x|}{\sqrt{2\pi s^3}} e^{-\frac{m_1^2 s}{2}-m_1x-\frac{x^2}{2s}}  ds dl dy\nonumber
%\\
%& \le  & \int_0^\infty e^{-\frac{m_*^2t}{2}-lp(m_1+m_2)+lm_2+m_1y-m_1x} \int_0^t  p^0(t-s, 0, y, l) \frac{|x|}{\sqrt{2\pi s^3}}e^{-\frac{x^2}{2s}}  ds  dl dy \nonumber
%\\
%& = & \int_0^\infty e^{-\frac{m_*^2t}{2}-lp(m_1+m_2)+lm_2+m_1y-m_1x}\int_0^t  p^0(t-s,0,y,l) \IP_x\left[\sigma^0_{\{0\}}\in ds\right] dl dy \nonumber
%\\
&\le & \int_0^\infty e^{-\frac{m_*^2t}{2}-lp(m_1+m_2)+lm_2+m_1y-m_1x} p^0(t,x,y,l)dldy, \label{eq3.14}
\end{eqnarray}
as well as 
\begin{eqnarray}
&&  \int_0^t  \IP^0[Y_{t-s}\in dy]\; \IP^x[\sigma_{\{0\}} \in ds]  \nonumber
\\
&\ge &  \int_0^\infty e^{-\frac{(m^*)^2t}{2}-lp(m_1+m_2)+lm_2+m_1y-m_1x} p^0(t,x,y,l)dldy, \label{eq3.15}
\end{eqnarray}
\cite[Corollary 3.3]{ABTWW, ABTWW1} gives
\begin{equation}\label{eq3.16}
p^0(t,x,y,l)=\frac{2p(l+y+x)}{\sqrt{2\pi t^3}} e^{-\frac{(l+y+x)^2}{2t}}, \quad l>0, y\le 0.
\end{equation}
Replacing the $p^0(t,x,y,l)$ in \eqref{eq3.14} and \eqref{eq3.15} with the right hand side of \eqref{eq3.16} yields
\begin{eqnarray} 
&&\frac{2p}{\sqrt{2\pi t^3}} e^{-\frac{(m^*)^2t}{2}+m_1(y-x)}\int_0^\infty (l+y+x)  e^{-lp(m_1+m_2)+lm_2-\frac{(l+y+x)^2}{2t}}dldy \nonumber 
\\
&\le &\int_0^t  \IP^0[Y_{t-s}\in dy]\; \IP^x[\sigma_{\{0\}} \in ds] \nonumber
\\
&\le & \frac{2p}{\sqrt{2\pi t^3}} e^{-\frac{m_*^2t}{2}+m_1(y-x)}\int_0^\infty (l+y+x)  e^{-lp(m_1+m_2)+lm_2-\frac{(l+y+x)^2}{2t}}dldy. \label{case2-eq3.14}
\end{eqnarray}
In view of  \eqref{case2-density},  \eqref{case2-eq3.14} combined with \eqref{density-of-q} gives the desired result for Case 2.

To deal with cases (iii) and (iv) in Theorem \ref{main-result}, i.e., when $x<0$, we denote by $\widetilde{Y}:= -Y$. By \cite[p.218, Theorem 7.1]{KS}, noticing that the symmetric semimartingale local time of $Y$ and $\wt{Y}$ (whose existence is guaranteed also by \cite[p.218, Theorem 7.1]{KS})  at $0$ are the same, we have 
\begin{eqnarray}
d\widetilde{Y}_t  &=& -m(Y_t)dt -dB_t-(2p-1)dL_t^0(Y) \nonumber
\\
&= &  -\left(m_1\cdot \mathbf{1}_{\{Y_t\ge 0\}} +m_2\cdot \mathbf{1}_{\{Y_t<0\}}\right)dt -dB_t -(2p-1)dL_t^0(\wt{Y}) \nonumber
\\
&=& d\wt{B}_t + \left( (-m_2)\cdot \mathbf{1}_{\{\wt{Y}_t>0\}} +(-m_1)\cdot \mathbf{1}_{\{\wt{Y}_t\le 0\}}  \right)dt +\left(2(1-p)-1\right) dL^0_t(\wt{Y}),
\end{eqnarray}
where $\wt{B}:=-B$ is also a standard Brownian motion.  Since $\wt{Y}$ also does not have sojourn time at $0$,  we can characterize $\wt{Y}$ by the following SDE:
\begin{equation}\label{sde-Y-ti}
d\widetilde{Y}_t =dB_t + \left( (-m_2)\cdot \mathbf{1}_{\{\wt{Y}_t\ge 0\}} +(-m_1)\cdot \mathbf{1}_{\{\wt{Y}_t<0\}}  \right)dt +\left(2(1-p)-1\right) dL^0_t(\wt{Y}).
\end{equation}
To complete cases  (iii) and (iv)  in Theorem \ref{main-result},  it suffices to notice that  when $x<0$,    $p^Y(t,x,y)  = p^{\wt{Y}}(t,-x,-y)$ can be found by replacing $m_1, m_2$ and $p$ with $-m_2, -m_1$, and $1-p$  in  cases (i) and (ii), respectively.

\vskip 0.3truein

\noindent {\bf Shuwen Lou}

\smallskip \noindent
Department of Mathematics and Statistics, Loyola University Chicago,
\noindent
Chicago, IL 60660, USA

\noindent
E-mail:  \texttt{slou1@luc.edu}

 \end{document}